\documentclass[a4paper,10pt]{article}
\usepackage{amsmath,amsfonts,amssymb,amsthm}
\usepackage{amsrefs}
\newcommand{\Rset}{\mathbb{R}}
\newcommand{\Nset}{\mathbb{N}}

\newcommand{\PP}{\mathbb{P}}
\newcommand{\EE}{\mathbb{E}}

\newtheorem{theorem}{Theorem}[section]
\newtheorem{lemma}[theorem]{Lemma}
\newtheorem{definition}[theorem]{Definition}
\newtheorem{proposition}[theorem]{Proposition}
\numberwithin{equation}{section}
\title{On the dynamic programming approach for the $3D$ Navier-Stokes equations}
\author{Luigi Manca \footnote{Scuola Normale Superiore, Piazza dei Cavalieri 7, 56126 Pisa, Italy. \emph{E-mail: l.manca@sns.it}}}

\begin{document} 

\maketitle

\begin{abstract}
The dynamic programming approach for the control of a $3D$ flow governed by the stochastic Navier-Stokes equations for incompressible fluid in a bounded domain is studied.
By a compactness argument, existence of solutions for the associated  Hamilton-Jacobi-Bellman equation is proved.
Finally, existence of an optimal control through the feedback formula and of an optimal state is discussed.\\
\bigskip

\noindent
\begin{center} \textbf{R\'esum\'e}  \end{center}
\bigskip

Nous \'etudions la programmation dynamique du contr\^ole d'un flux tridimensionnel gouvern\'e par les \'equations stochastiques de Navier-Stokes qui concernent un fluide incompressible dans un domaine born\'e.
Nous d\'emontrons l'existence de solutions pour l'\'equation associ\'ee de Hamilton-Jacobi-Bellman  par un argument de compacticit\'e.
Enfin nous examinons l'existence d'un contr\^ole optimal et d'un \'etat optimal au moyen de la formule de feedback.\bigskip

\noindent
{\em MSC:} 76D05; 76D55; 49L20\\
{\em Keywords:} Navier-Stokes equations; dynamic programming; Hamilton-Jacobi-Bellman equations
\end{abstract}

\section{Introduction}
In this article we study the dynamic programming approach for the control of a three dimensional turbulent flow
governed by the stochastic Navier-Stokes equations for incompressible fluids.
The unknows are the velocity field $U(\xi,t) = (U_1(\xi,t), U_2(\xi,t), U_3(\xi,t))$ and the pressure $p(\xi,t)$, where $t\in [0,T]$ and $\xi \in D$, with $D\subset \Rset^3$ open and bounded. 
$U(\xi,t)$, $p(\xi,t)$ satisfy the equation
\begin{equation*}
\begin{cases}
  \frac{\partial U}{\partial t} + (U\cdot\nabla)U + \nabla p = 
      \nu \Delta U+z+\dot{\eta}\quad\text{in  } D \\ 
  \text{div} U = 0, \quad\text{in  } D \\
  U=0 \quad\text{in  } \partial D \\
  U|_{t=0} = u_0 \quad\text{in  } D.
\end{cases}
\end{equation*}
where the control $z=z(\xi,t)$ is a bounded random variable, $\nu$ is the kinematic viscosity and $\eta$ is a random variable of white noise type.
Since $\nu$ does not play a particular role, with no loss of generality we can assume that $\nu = 1$.
We look for a solution with values in the Hilbert space $H$ of the square integrable and divergence free functions $f:D\to \Rset^3$.  

We consider a cost functional of the form
$$
  J(z)=\EE \bigg[ \int_D\int_0^T\big(\Phi( U(\xi,t)) + \frac{1}{2}|z(\xi,t)|^2\big)dt + \varphi(U(\xi,T))\big)d\xi\bigg],
$$
where  $T>0$ and $\Phi,\varphi:\Rset^3\to \Rset^+$ are given functions.

The idea is that $J(z)$ measures the amount of turbulence of the system.
So, in \cites{DPD00, DPD001} it is proposed to take $\varphi$ bounded and $\Phi(x)= |\nabla \times x|^2$, where $\nabla\times x$ is the rotational of the three dimensional fields $x$.
In this paper the running cost $\Phi$ satisfies stronger conditions, but we think that our assumptions will be not very restrictive.

Many articles are been devoted to this control problem, and it is proved that the optimal control value is a viscosity solution of the Hamilton-Jacobi-Bellman(HJB) equation associated to the problem (see, for instance, \cites{DPD00, DPD001, GSS05} and references therein).
Unfortunately, viscosity solutions are not smooth enough to fully justify the dynamic programming approach.

We follow a strategy proposed in \cites{CDP90, CDP92} \ to get smooth solutions.
In \cites{DPD00, DPD001} this strategy has been succesfully implemented to study the dynamic programming approach for the Burgers equations and for the $2D$ Navier-Stokes equations. 

After  delicate a priori estimates on the Galerkin approximated problem, we are able to find a solution for the HJB equation by compactness argument. 
Unfortunately, as for the uncontrolled equation, we are not able to prove uniqueness.
Moreover, due to the lack of informations on the differentiability of the flow with respect to the space variable, we are not able to apply verification theorems (cf. \cite{G02}).
Then we shall justify the dynamic programming approach only  for determinated classes of solution of the controlled equation, which depend by a given solution of the HJB equation.

\section{Notations}
Let $D\subset \Rset^3$ be a bounded open set with regular boundary $\partial D$ in $\Rset^3$ and let  $L^2(D)$ be the set of the real valued square Lebesgue integrable functions on $D$. We denote by $H^k(D)$ the usual Sobolev spaces, and by $H_0^1(D)$ the space of all functions in $H^1(D)$ with vanish on the boundary $\partial D$.
We introduce the Hilbert spaces
$$
      H=\left\{x\in (L^2(D))^3:
         \textrm{div}  x = 0 \text{ in } D,\, x = 0\text{ in } \partial D \right\},
$$
$$
     V=\left\{x\in (H_0^1(D))^3:\textrm{div} x = 0 \text{ in } D \right\},
$$
where $n$ denotes the normal unit vector on $\partial D$. 
$H$ (resp. $V$) is endowed with the inner product and norm of $(L^2(D))^3$ (resp. $(H_0^1(D))^3$) denoted by $(\cdot,\cdot)$ and $|\cdot|$ (resp. $((\cdot,\cdot))$ and $\|\cdot\|$).
Moreover, we introduce the unbounded self-adjoint operator
$$
                   A=P\Delta,\quad D(A)= (H^2(D))^3\cap (H^1_0(D))^3\cap H,
$$
where $P$ is the orthogonal projector of $(L^2(D))^3$ onto $H$ and  the operator $b$ is defined by
$$
              b(x,y)= P\big((x\cdot \nabla)y\big),\quad b(x)=b(x,x)\quad x,y\in V.
$$
$W$ is a cylindrical Wiener process defined on a stochastic basis $(\Omega,\mathcal{F},\mathcal{F}_{t\geq 0},\PP)$ with values in $H$.
The operator $Q$ is symmetric, nonnegative, of trace class and such that $\ker Q=\{0\}$.
The control $z$ is chosen in the space of adapted processes
$$
        \mathcal{M}_R = \{z \in L_W^2(\Omega\times[0,T];H), |z|\leq R \},
$$
for a fixed $R>0$, and it is subject to a linear operator $B\in \mathcal{L}(H)$
which will be specified below.
We study our control problem with initial value in $D(A)$.
This choise will be clearify in the following.

Thanks to the introduced notations, we can write the problem in the abstract form\footnote{When it will be necessary to emphasize the dependence of $X$ by the initial condition $x$, we write $X(t,x)$ instead of $X(t)$.}
\begin{equation} \label{contr}
\begin{cases} dX(t)= (AX(t) + b(X(t)) + Bz(t))dt + Q^{1/2}dW(t), \\
              X(0)=x,\, x\in D(A).
\end{cases}
\end{equation}
Equation \eqref{contr} is associated to the {\em cost function}
\begin{equation*} 
   J(X,z) = 
  \EE\bigg[ \int_0^T \big(\Phi(X(t)) + \frac{1}{2}|z(t)|^2\big)ds + \varphi(X(T))\bigg],
\end{equation*}
where $T>0$ is  fixed and $\Phi,\varphi:D(A)\to\Rset$ are suitable nonnegative functions.
We stress that the cost function $J$ depends also by the solution of \eqref{contr} since we have no informations about its uniqueness. 
It is easy to see, by the It\^o formula, that since the covariance operator $Q$ is of trace class then $J(X,z)$ is well defined.
We aim at finding an {\em optimal control} $z^* \in \mathcal{M}_R$ and an {\em optimal state} $X^*$ which minimize $J$:
$$
     J(X^*,z^*)= \min\big\{J(X,z):z\in\mathcal{M}_R,\,X\text{ is a solution of \eqref{contr}}\big\}.
$$
We follow the dynamic programming approach to solve this problem.
Let the Hamiltonian $F$ be defined on $H$ by
\begin{equation*} 
 F(p)= \begin{cases} 
                     \frac{1}{2}|p|^2,    &\text{ if } |p|\leq R, \\
                     |p|R-\frac{1}{2}R^2, &\text{ if } |p|> R.
       \end{cases}
\end{equation*}
The Hamilton-Jacobi-Bellman equation associated with our control problem is
\begin{equation} \label{HJB}
\begin{cases}
    D_t u= 
     \frac{1}{2}\textrm{Tr}[Qu_{xx}]+(Ax+b(x),u_x)-F(B^*u_x) + \Phi(x),\\ 
    u(0,x) = \varphi(x),\,x\in D(A),
\end{cases}
\end{equation}
where the subscript $x$ means the differential with respect to $x$. 
If we are able to find a smooth solution $u$ of \eqref{HJB} then the optimal control is given by the \emph{feedback formula} (cf. \cite{DPZ02})
\begin{equation} \label{z^*}
     z^*(t) = -D_pF(B^*u_x(T-t,\widetilde{X}(t))),
\end{equation}
where the optimal state $\widetilde{X}$ is the solution of the \emph{closed loop} equation
\begin{equation} \label{loop}
\begin{cases} d\widetilde{X}(t) 
            = (A\widetilde{X}(t) + b(\widetilde{X}(t)) -D_pF(B^*u_x(T-t,\widetilde{X}(t)))dt + Q^{1/2}dW(t), \\
              \widetilde{X}(0)=x \in D(A).
\end{cases}
\end{equation}
Due to the lack of uniqueness of $u(t,x)$ we shall solve the control problem only for some classes of solution of \eqref{contr}.
This discussion is detailed in the last section.
We shall solve equation \eqref{HJB} under the mild form
\begin{eqnarray} 
 u(t,x) &=& R_t\varphi(x) + \int_0^tR_{t-s}(b(\cdot),u_x(s,\cdot))(x)ds \notag \\ \label{R_t}
         &&- \int_0^t R_{t-s}F(B^*u_x(s,\cdot))(x)ds + \int_0^t R_{t-s}\Phi(x)ds,
 \end{eqnarray}
where $R_t$ is the Ornstein-Uhlembeck semigroup defined by
$$
   R_t\varphi(x) = \EE[\varphi(Z(t,x))],\quad \varphi \in C_b(D(A);\Rset)
$$
and $Z(t,x)$ is the solution of the stochastic equation
$$
\begin{cases}
      dZ(t) = AZ(t)dt + Q^{1/2}dW(t), &t>0 ,\\
      Z(0)=x,                         &x\in D(A).
\end{cases}
$$
%
%
%
%
\section{Galerkin approximations}
We introduce the usual Galerkin approximations of equation \eqref{contr}. 
For $m\in \Nset$, we define by $P_m$ the projector of $H$ onto the space spanned by the first $m$ eigenvectors of $A$.
Then we set, for $x\in H$, 
$b_m(x) = P_m b(P_m x)$, $Q_m=P_m Q P_m$ and $B_m=P_mBP_m$.
We consider the control problem in finite dimension for the approximated equation
\begin{equation} \label{contrm}
\begin{cases} dX^m(t) 
            = \big(AX^m(t) + b_m(X^m(t))+ B_m z(t)\big)dt + Q^{1/2}_mdW(t), \\
              X^m(0)=P_m x,\, x \in D(A),
\end{cases}
\end{equation}
which consists in minimizing the cost function
\begin{equation*} 
 J_m(z_m) = 
  \EE\bigg[ \int_0^T \big(\Phi(X^m(t)) + 
  \frac{1}{2}|P_mz(t)|^2\big)ds + \varphi(X^m(T))\bigg],
\end{equation*}
where $z \in \mathcal{M}_R$. .
Then, we consider the approximated equation of 
\eqref{R_t}, i.e.
\begin{eqnarray} 
 u^m(t,x) =&& R_t^m\varphi(x) + \int_0^tR_{t-s}^m(b_m(\cdot),u_x^m(s,\cdot))(x)ds \notag\\
 \label{R_tm}
         &&- \int_0^t R_{t-s}^mF(B_m^*u_x^m(s,\cdot))(x)ds + \int_0^t R_{t-s}^m\Phi(x)ds,
 \end{eqnarray}
where $\{R_t^m\}_{t\geq0}$ is the Galerkin approximation of the Ornstein-Uhlembeck semigroup.
In order to approximate \eqref{R_t} by \eqref{R_tm}, we need some a priori bounds for $u^m(t,x)$.
So, we consider the Feynman-Kac semigroup $\{S_t^m\}_{t\geq 0}$ defined by
\begin{equation} \label{FKsgr}
   S_t^m\varphi(x) = \EE\Big[\text{\rm e}^{-K\int_0^t|AY^m(s,x)|^2ds}\varphi(Y^m(t,x))\Big],
\end{equation}
where $Y_m(t,x)$ is the solution of the finite dimensional stochastic equation
\begin{equation*} 
\begin{cases}
       dY^m(t)  = (AY^m(t) + b_m(Y^m(t)))dt + Q^{1/2}_mdW(t), \\
         Y^m(0) = P_m x,\quad x \in D(A).
\end{cases}
\end{equation*}
\eqref{R_tm} has an unique solution $u^m(t,x)$, which is also the solution of
\begin{eqnarray} \label{eq4m}
 u^m(t,x)&=&S_t^m\varphi(x) + \int_0^t S_{t-s}^m(|A\cdot|^2 u^m(s,\cdot))(x)ds \notag \\
 &&-\int_0^t S_{t-s}^m F_{}(B_m^*u_x^m(s,\cdot))(x)ds + \int_0^t S_{t-s}^m \Phi(x)ds,
\end{eqnarray}
and of
\begin{equation} \label{HJBm}
\begin{cases}
    D_t u^m(t,x)
     =\frac{1}{2}\textrm{Tr}[Q_mu_{xx}^m(t,x)]+(Ax+b_m(x),u_x^m(t,x))\\
        \qquad\qquad\qquad-F(B_m^*u_x^m(t,x)) + \Phi(x),\\
    u^m(0,x) = \varphi(x),
\end{cases}
\end{equation}
where $(t,x)\in [0,T]\times P_mH$.
Thanks to this property, we shall be able to find the a priori bounds that we need. 
By a classical computation based on the It\^o formula (cf. \cite{DPZ02})
 we find that 
the optimal control $z_m^*$ is obtained by taking
$$
  z_m^*(t) = -D_pF_{}(B_m^*u_x^m(T-t,\widetilde{X}^m(t))),
$$
where $\widetilde{X}^m(t,x)$ is the solution of the closed loop equation
\begin{equation} \label{loopm}
\begin{cases} d\widetilde{X}^m(t) 
      = \big(A\widetilde{X}^m(t) + b_m(\widetilde{X}^m(t)) \\ \qquad\qquad\quad-D_pF_{}(B_m^*u_x^m(T-t,\widetilde{X}^m(t)))\big)dt + Q^{1/2}_mdW(t), \\
       \widetilde{X}^m(0)=P_mx,\, x\in D(A).
\end{cases}
\end{equation}
Since \eqref{loopm} is in finite dimension, it is easy to prove existence and uniqueness of a solution. 

We recall that a detailed survey on HJB equations on finite dimensional spaces may be found in \cite{FS}. 
The infinite dimensional problem has been recently developed, see \cite{DPZ02}.
%
%
%
\section{Functional spaces}
Let  $(E,|\cdot|_E )$ be a Banach space. 
If $\varphi:D(A) \to E$ and $h\in H$ we set
$$
     D\varphi(x)\cdot h = \lim_{s \to 0} \frac{\varphi(x +sh) - \varphi(x)}{s},
$$
when the limit in $E$ exists. 
The space $D(A)$ is endowed with the graph norm.
We define the following functional spaces:
\begin{itemize}
\item $C_b(D(A);E)$ is the space of all continuous and bounded functions from $D(A)$ to $E$, endowed with the norm
$$
   |\varphi |_0 := \sup_{x\in D(A)}|\varphi(x)|_E, \quad \varphi \in C_b(D(A);E).
$$
\item For any $k\in \Nset$, $C_k(D(A);E)$ is the space of all continuous functions from $D(A)$ to $E$ such that
$$
   |\varphi|_{k,A}:= \sup_{x\in D(A)}\frac{|\varphi(x)|_E}{(1+|Ax|)^k} <\infty.
$$
\item For any $k\in \Nset$, $C_k^1(D(A);E)$ is the space of all functions of $C_k(D(A);E)$ such that
$$
 |\varphi|_{k,A,1}:= 
          \sup_{x\in D(A)}\frac{|(-A)^{-1}D\varphi(x)|_E}{(1+|Ax|)^k} <\infty.
$$
\item  For any $k\in \Nset$, $\mathcal{E}_k(D(A);E)$ is the space of all function $\varphi$ $\in C_k(D(A);E)$ such that
$$
   |\varphi|_{\mathcal{E}_k} = \sup_{\substack{ x,y\in D(A)\\x\not= y}}
     \frac{|\varphi(x) - \varphi(y)|}{|A(x-y)|(1+|Ax|+|Ay|)^k} <\infty.                                  
$$
\item For any $k\in \Nset$, $\alpha\in (0,1)$ the set $C_{k,(-A)^\alpha}^1(D(A);E)$ denotes  the space of all functions of $C_k^1(D(A);E)$ such that
$$
 |\varphi|_{k,(-A)^\alpha,1}:= 
          \sup_{x\in D(A)}\frac{|(-A)^{-\alpha}D\varphi(x)|_E}{(1+|Ax|)^k} <\infty.
$$
\end{itemize}
For any $\delta>0$ we denote by  $D((-A)^{-\delta})$ the dual space of $D((-A)^\delta)$. 
We shall identify the space $H$ with its dual $H^*$.
Hence, the embeddings $D((-A)^\delta)$ $\subset$ $H$ $\subset$ $D((-A)^{-\delta})$ hold.
%
%
%
\section{Hypothesis on the operators $Q$, $B$}
Following \cites{DPD03, DO} we assume that 
\begin{equation} \label{Q1}
 Tr[(-A)^{1+g}Q] < +\infty
\end{equation}
for some $g>0$ and that 
\begin{equation} \label{Q2}
 |Q^{-1/2}x| \leq c_r|(-A)^rx|,\quad \forall x \in D((-A)^r),
\end{equation}
for some $r\in (1,3/2)$ and $c_r >0$.
Then, we shall denote by $\varepsilon$ the quantity
\begin{equation} \label{varepsilon}
 \varepsilon = \frac{3-2r}{2}.
\end{equation}
Notice that by the hypothesis on $r$ we have $\varepsilon \in (0,1/2)$.
If $Z(t)$ denotes the stochastic convolution with covariance operator $Q$, i.e. $Z(t)$ is the solution of the linear stochastic equation in $H$
$$
 Z(t)=\int_0^t\text{\rm e}^{(t-s)A}\sqrt{Q}dW(s), \quad t\geq 0,
$$
then hypothesis \eqref{Q1}  implies that
$$
   \EE\big[|(-A)^{1+g/2}Z(t)|^2\big] = \textrm{Tr}\big[(-A)^{1+g}Q\big] < \infty.
$$
Since $Z(t)$ is gaussian, by the factorization method (see, for instance, \cite{DPZ}) 
it follows that for any $p>0$ it holds
$$
  \EE\bigg[\sup_{t\in[0,T]}|(-A)^{1+g/2}Z(t)|^p\bigg] \leq c(p)\textrm{Tr}\big[(-A)^{1+g}Q\big]^\frac{p}{2}.
$$
Consequently, setting 
$$
   Z_m(t) = \int_0^t\text{\rm e}^{(t-s)A}\sqrt{Q}_mdW(s)
$$
it is clear that 
\begin{equation*} 
  \EE\bigg[\sup_{t\in[0,T]}|(-A)^{1+g/2}Z_m(t)|^p\bigg] \leq c(p)\textrm{Tr}\big[(-A)^{1+g}Q\big]^\frac{p}{2}.
\end{equation*}
On the linear operator $B$ we suppose that
$$
 B: H\longrightarrow D((-A)^{\gamma}) 
$$
for some $\gamma >1-\varepsilon$. 
This implies that for some $c> 0$ and for all $z\in D((-A)^{-\gamma})$ it holds
\begin{equation} \label{B}
 |B^*z|\leq c|(-A)^{-\gamma}z|. 
\end{equation}
%
%
%
%
\section{A priori estimates I}
In this section we prove some useful estimates on the Feynman-Kac semigroup \eqref{FKsgr}.
We  omit the proof of well-known results.
We use the following estimates on the bilinear operator $b(x,y)$ (see, for instance, \cites{DPD03, Tem77}).
\begin{lemma} \label{lemma2.1}
There exists $c>0$ such that for all $x,y \in D(A)$, $z\in H$.
$$
(b(x,y),(-A)^{1/2}z)\leq c|Ax||Ay||z|.
$$
\end{lemma}
Detailed proof of the following lemmata may be found in \cite{DPD03}.
\begin{lemma} \label{lemmaa}
Let $k\in \Nset$, $\varphi \in C_k(D(A);\Rset)$ $\cap$ $C_{k,(-A)^{1/2}}(D(A);\Rset)$ and $\alpha\in (1/2,1)$. 
Then there exists $c>0$ such that if $K$ is sufficiently large we have
 \begin{eqnarray*} 
  |S^m_t\varphi|_{k,(-A)^{1/2},1} &\leq& ct^{-1/2}(|\varphi|_{k,A}+|\varphi|_{k,A,1}),\quad t>0, \\
  |S^m_t\varphi|_{k,(-A)^{\alpha},1} &\leq& ct^{\alpha-1}(|\varphi|_{k,A}+|\varphi|_{k,A,1}) \quad t>0.
 \end{eqnarray*}
 \end{lemma}
 \begin{proof}
 The first estimate is proved in Lemma 4.4 in \cite{DPD03}.
 The second one follows by interpolation between the first one and the estimate
 $$
   |S_t^m|_{k,A,1}\leq c(|\varphi|_{k,A}+ |\varphi|_{k,A,1}),
 $$
 which is proved in Lemma 4.2 of \cite{DPD03}.
 \end{proof}
\begin{lemma} \label{lemmab}
Let $k\in \Nset$, $\varphi \in C_k(D(A);\Rset)$, $\alpha\in (1/2,1)$.
Then there exists $c>0$ such that if $K$ is sufficiently large we have
\begin{equation*} 
 |S^m_t\varphi|_{k,(-A)^\alpha,1} \leq c_1 (1+t^{\alpha+\varepsilon -2})|\varphi|_{k,A},\quad t>0,
\end{equation*}
where $\varepsilon$ is defined by \eqref{varepsilon}.
\end{lemma}
\begin{proof}
By Lemma 4.1 of \cite{DPD03},  we know that there exists $c>0$ such that 
\begin{equation*} 
  |S^m_t\varphi|_{k,A,1} \leq c(1+t^{\varepsilon -1})|\varphi|_{k,A},\quad t>0,\,m\in\Nset.
 \end{equation*}
Hence, taking into account Lemma \ref{lemmaa} and the semigroup property of $S_t^m$, we have 
\begin{eqnarray*}
  &&|S^m_t\varphi|_{k,(-A)^\alpha,1} \\
  &&\quad= |S^m_{t/2}S^m_{t/2}\varphi|_{k,(-A)^\alpha,1} 
     \leq c(1+(t/2)^{\alpha-1})(|S^m_{t/2}\varphi|_{k,A}+|S^m_{t/2}\varphi|_{k,A,1})\\
  &&\quad   \leq c(1+(t/2)^{\alpha-1})(|\varphi|_{k,A}+(1+(t/2)^{\varepsilon-1})|\varphi|_{k,A})
\end{eqnarray*}
which implies the result.
\end{proof}
\begin{lemma} \label{lemmac}
 Let $k\in\Nset$, $\varphi \in C^1_k(D(A);\Rset)$, $\alpha\in(1/2,1)$ and  $\sigma \in (3/4,1)$.  
 Then there exists $c>0$ such that for all $m\in\Nset$, $t> 0$  it holds
\begin{equation*}
 |(-A)^{-\alpha}D^2S_t^m\varphi(\cdot)(-A)^{-\sigma}|_{k,A} \leq c(1+t^{\sigma+\alpha+\varepsilon-3})(|\varphi|_{k,A}+|\varphi|_{k,A,1}).
\end{equation*}
\end{lemma}
\begin{proof}
The result follows by  Lemma 4.5 and Lemma 4.8 in \cite{DPD03} and by arguing as in the previous Lemma.
\end{proof}
%
%
\section{A priori estimates II} 
In this section, we assume that $\varphi$, $\Phi$, $\alpha$ satisfy the following conditions
\begin{equation}\label{varphi}
\begin{array}{l}
  \varphi(x),\,\Phi(x)\geq0 \quad x\in D(A);\\{}\\
  \varphi,\,\Phi\in  C_b(D(A);\Rset)\cap\mathcal{E}_{2}(D(A);\Rset) \\{}
  \\
  \alpha\in(1-\varepsilon,1),\quad \alpha \leq \gamma
\end{array}
\end{equation}   
 where  $p\in\Nset$ is fixed, $c>0$ and $\varepsilon$, $\gamma$ are  defined in \eqref{varepsilon}, \eqref{B} respectively. 
We are going to establish some estimates on $u^m(t,x)$, $Du^m(t,x)$, $D^2u^m(t,x)$. 
\begin{proposition}\label{E.0}
   There exists $c >0$ such that for all $m\in\Nset$ it holds
  \begin{equation*} 
  \sup_{t\in[0,T]}|u^m(t,\cdot)|_0 \leq c(|\varphi|_0+|\Phi|_0)
 \end{equation*}
\end{proposition}
\begin{proof}
Since $u^m$ is the function associated with the approximated control problem we deduce that
 \begin{eqnarray*}
    0&\leq& u^m(T,x) = \min_{z\in\mathcal{M}_R\cap L^2_W(\Omega \times [0,T];P_m H)} J_m(z)  
    \leq J_m(0) \\
   &=& c\EE\left[\int_0^T \Phi(X^m(s,x))ds +\varphi(X^m(T,x))\right]\\
    &\leq& T|\Phi|_0+|\varphi|_0. 
 \end{eqnarray*}
  We could consider the same control problem in $[0,t]$ and obtain 
 $$
   0\leq u^m(t,x) \leq |\varphi|_0+T|\Phi|_0
 $$
 This completes the proof.
\end{proof}
 \begin{proposition} \label{E.1}
There exists $c>0$ such that for all $m\in \Nset$ it holds
\begin{equation*} 
  |u^m(t,\cdot)|_{2,(-A)^\alpha,1} \leq ct^{\alpha-1}(|\varphi|_0 +|\varphi|_{\mathcal{E}_{2}}+|\Phi|_0+|\Phi|_{\mathcal{E}_{2}}).
\end{equation*}
\end{proposition}
\begin{proof}
Let us assume that $\varphi,\Phi \in C_b(D(A);\Rset)\cap C_{2}^1(D(A);\Rset)$.
Taking into account \eqref{eq4m} we have
$$
  |u^m(t,\cdot)|_{2,(-A)^\alpha,1} \leq I_1+I_2+I_3+I_4
 $$
 where
\begin{eqnarray*}
   I_1 &=& |S_t^m\varphi(\cdot)|_{2,(-A)^\alpha,1},\\
   I_2 &=& \int_0^t |S_{t-s}^m(|A\cdot|^2 u^m(s,\cdot))|_{2,(-A)^\alpha,1}ds, \\
   I_3 &=& \int_0^t |S_{t-s}^mF(B_m^*u_x^m(s,\cdot))|_{2,(-A)^\alpha,1}ds, \\
   I_4 &=& \int_0^t |S_{t-s}^m\Phi(\cdot)|_{2,(-A)^\alpha,1}ds.
\end{eqnarray*}
 $I_1$ and $I_4$ are estimated as in Lemma \ref{lemmaa}. 
For $I_2$ we have that
$$
   ||A\cdot|^2 u^m(s,\cdot)|_{2,A} \leq c| u^m(s,\cdot)|_0 \leq c(|\varphi|_0+|\Phi|_0)
$$
by Proposition \ref{E.0}.
Hence, taking into account Lemma \ref{lemmab} and that $\alpha+\varepsilon>1$, it holds
\begin{equation*} 
 I_2 \leq  c \int_0^t(1+(t-s)^{\alpha+\varepsilon -2})||A\cdot|^2 u^m(s,\cdot)|_{2,A}ds \leq c(|\varphi|_0+|\Phi|_0).
\end{equation*}
For $I_3$ we have that 
 $$
  |S_t^m F_{}(B_m^*u_x^m(s,\cdot))|_{2,(-A)^\alpha,1} \leq c
   (1+t^{\alpha+\varepsilon -2})|F_{}(B_m^*u_x^m(s,\cdot))|_{2,A}.
 $$
 by Lemma \ref{lemmab}.
 Since $\alpha \leq \gamma$, by \eqref{B} it follows
 \begin{eqnarray*}
  &&|F(B_m^*u_x^m(s,\cdot))|_{2,A} \leq   c|(-A)^{-\gamma}u_x^m(s,\cdot)|_{2,A}\\
  &&\quad\leq c|(-A)^{-\alpha} u_x^m(s,\cdot)|_{2,A} = c|u^m(s,\cdot)|_{2,(-A)^\alpha,1},
 \end{eqnarray*}
 and so we find
\begin{equation*} 
 I_3\leq c\int_0^t(1+(t-s)^{\alpha+\varepsilon -2})
 |u^m(s,\cdot)|_{2,(-A)^\alpha,1}ds.
\end{equation*}
Finally, by gathering all the estimates  on $I_1$, $I_2$, $I_3$, $I_4$ 
   we find that for some $c>0$ it holds
\begin{eqnarray*}
   |u^m(t,\cdot)|_{2,(-A)^\alpha,1}  &\leq& ct^{\alpha-1}(|\varphi|_0+|\varphi|_{2,A,1}+|\Phi|_0+|\Phi|_{2,A,1})  \\
  &&+   c \int_0^t (1+(t-s)^{\alpha+\varepsilon -2}) 
            |u^m(s,\cdot)|_{2,(-A)^\alpha,1}ds.
\end{eqnarray*}
 Then,  by  Gronwall's lemma (see, for instance,  \cite{Hen}),  we have
 \begin{equation*}
   |u^m(t,\cdot)|_{2,(-A)^\alpha,1} 
 \leq ct^{\alpha-1}(|\varphi|_0+|\varphi|_{2,A,1}+|\Phi|_0+|\Phi|_{2,A,1})
\end{equation*}
 since $\alpha+\varepsilon >1$.
 Now notice that all the estimates above are done in the finite dimensional space $P_mH$.
Hence, if $\varphi$, $\Phi$ $\in C_b(D(A);\Rset)\cap \mathcal{E}_{2}(D(A);\Rset)$, we obtain the result by approximating uniformly $\varphi$, $\Phi$ by functions in   $C_b(D(A);\Rset)$ $\cap$ $C_{2}^1(D(A);\Rset)$. 
\end{proof}
The following two results will be proved with a similar argument, by using Lemma \ref{lemmac}.
\begin{proposition} \label{E.2}
  Let $\sigma\in(3/4,1)$ such that $\sigma+\alpha+\varepsilon >2$, where $\alpha$ and $\varepsilon$ are defined in \eqref{varphi}, \eqref{varepsilon} respectively.
  Then there exists $c>0$ such that for all  $m\in \Nset$ it holds
\begin{eqnarray*} 
  &&|(-A)^{-\alpha}D^2u^m(t,\cdot)(-A)^{-\sigma}|_{4,A}    \\
  &&\quad  \leq ct^{\sigma+\alpha+\varepsilon-3}
     (|\varphi|_{0}+|\varphi|_{\mathcal{E}_{2}}+|\Phi|_{0}+|\Phi|_{\mathcal{E}_{2}}). 
\end{eqnarray*}
\end{proposition}
\begin{proof}
Notice that by the approximation argument described in Proposition \ref{E.1} it is sufficient to prove the claim for $\varphi,\Phi$ $\in$ $C_b(D(A);\Rset)$ $\cap$ $C_{2}^1(D(A);\Rset)$.
By \eqref{eq4m} we write
$$
|(-A)^{-\alpha}D^2u^m(t,\cdot)(-A)^{-\sigma}|_{4,A} \leq \sum_{i=1}^{4} J_i,
$$
 where
\begin{eqnarray*}
   J_1 &=& |(-A)^{-\alpha}D^2 S_t^m\varphi(\cdot)(-A)^{-\sigma}|_{4,A},\\
   J_2 &=& \int_0^t |(-A)^{-\alpha}D^2S_{t-s}^m(|A\cdot|^2 u^m(s,\cdot))(-A)^{-\sigma}|_{4,A}ds, \\
   J_3 &=& \int_0^t|(-A)^{-\alpha}D^2S_{t-s}^m F_{}(B_m^*u_x^m(s,\cdot))(-A)^{-\sigma}|_{4,A}ds, \\
   J_4 &=& \int_0^t|(-A)^{-\alpha}D^2 S_{t-s}^m \Phi(\cdot)(-A)^{-\sigma}|_{4,A}ds.
\end{eqnarray*}
$J_1$ and $J_4$ are estimated by Lemma  \ref{lemmac}       
For $J_2$ we have, by Lemma \ref{lemmac} 
\begin{eqnarray*}
 J_2 \leq&& c\int_0^t (t-s)^{\sigma+\alpha+\varepsilon-3}(||A\cdot|^2 u^m(s,\cdot)|_{4,A}+  ||A\cdot|^2 u^m(s,\cdot)|_{4,A,1})ds  \\
  \leq&& c\int_0^t (t-s)^{\sigma+\alpha+\varepsilon-3}(| u^m(s,\cdot)|_{2,A} + | u^m(s,\cdot)|_{3,A}+| u^m(s,\cdot)|_{2,A,1})ds.  
\end{eqnarray*}
Consequently,  taking into account Proposition \ref{E.0} 
and Proposition \ref{E.1} it follows 
\begin{eqnarray*} 
J_2\leq&& c(|\varphi|_{0}+|\varphi|_{2,A,1}+|\Phi|_{0}+|\Phi|_{2,A,1}) 
    \int_0^t (t-s)^{\sigma+\alpha+\varepsilon-3}s^{\alpha-1}ds \\
\leq&&c(|\varphi|_{0}+|\varphi|_{2,A,1}+|\Phi|_{0}+|\Phi|_{2,A,1}))t^{2\alpha+\sigma+\varepsilon-3}. 
\end{eqnarray*}
For $J_3$ we have, by Lemma \ref{lemmac}
\begin{equation*}
  J_3 \leq c\int_0^t (t-s)^{\alpha+\sigma+\varepsilon-3}(|F_{}(B_m^*u_x^m(s,\cdot))|_{4,A}+
      |F_{}(B_m^*u_x^m(s,\cdot))|_{4,A,1})ds.
\end{equation*}
Hence, since $|F_{}(B_m^*\varphi_x(\cdot))|\leq c|B^*\varphi_x(\cdot)|$, for $\alpha$ defined as in \eqref{varphi} we have that $J_3$ is bounded by
\begin{equation*}
 c\int_0^t (t-s)^{\alpha+\sigma+\varepsilon-3}(|u^m(s,\cdot)|_{4,(-A)^{\alpha},1} +|(-A)^{-\alpha} D^2u^m(t,\cdot)(-A)^{-1}|_{4,A})ds.
\end{equation*}
Taking into account Proposition \ref{E.1}, Proposition \ref{E.2} we find
\begin{eqnarray*} 
  J_3  &\leq& c(|\varphi|_{0}+|\varphi|_{2,A,1}+|\Phi|_{0}+|\Phi|_{2,A,1}))\\
  &&+c\int_0^t (t-s)^{\alpha+\sigma+\varepsilon-3}|(-A)^{-\alpha} D^2u^m(t,\cdot)(-A)^{-\sigma}|_{4,A})ds.
\end{eqnarray*}
So, the result follows by gathering the estimates on $J_1$, $J_2$, $J_3$, $J_4$ and by applying Gronwall's lemma.  
\end{proof}
\begin{proposition}\label{E.3}
There exists $c>0$ such that for all $m\in \Nset$ and $t>0$ it holds
\begin{equation*} 
|u^m(t,\cdot)|_{4,(-A)^{1/2},1} \leq ct^{-1/2}(|\varphi|_{0}+|\varphi|_{\mathcal{E}_{2}}+|\Phi|_{0}+|\Phi|_{\mathcal{E}_{2}}). 
\end{equation*}
\end{proposition}
\begin{proof}
By \eqref{eq4m} we write $|u^m(t,\cdot)|_{4,(-A)^{1/2},1}|\leq I_1+I_2+I_3+I_4$, where
\begin{eqnarray*}
   K_1 &=& |S_t^m\varphi(\cdot)|_{4,(-A)^{1/2},1},\\
   K_2 &=& \int_0^t |S_{t-s}^m(|A\cdot|^2 u(s,\cdot))|_{4,(-A)^{1/2},1}ds, \\
   K_3 &=& \int_0^t |S_{t-s}^mF_{}(B_m^*u_x(s,\cdot))|_{4,(-A)^{1/2},1}ds, \\
   K_4 &=& \int_0^t |S_{t-s}^m\Phi(\cdot)|_{4,(-A)^{1/2},1}ds.
\end{eqnarray*}
$K_1$, $K_2$ are estimated by Lemma \ref{lemmaa}.
For $K_2$ we have by Lemma \ref{lemmaa} 
\begin{eqnarray*}
 K_2 &\leq& c\int_0^t(t-s)^{-1/2}(||A\cdot|^2 u^m(s,\cdot))|_{4,A}+||A\cdot|^2 u^m(s,\cdot))|_{4,A,1}) ds\\
   &\leq& c\int_0^t(t-s)^{-1/2}(|u^m(s,\cdot)|_{2,A}+|u^m(s,\cdot)|_{3,A}+|u^m(s,\cdot)|_{2,A,1})ds.
\end{eqnarray*}
Then, taking into account Proposition \ref{E.2} and Proposition \ref{E.1} it follows
\begin{eqnarray*}
  K_2 &\leq& c(|\varphi|_{0}+|\varphi|_{2,A,1}+|\Phi|_{0}+|\Phi|_{2,A,1})\int_0^t (t-s)^{-1/2}s^{\alpha-1}ds\\
   &\leq& c(|\varphi|_{0}+|\varphi|_{2,A,1}+|\Phi|_{0}+|\Phi|_{2,A,1}) 
\end{eqnarray*}
since $\alpha>1/2$.
For $K_3$ we have
$$
 K_3\leq c\int_0^t(t-s)^{-1/2}(|F_{}(B_m^*u_x^m(t,\cdot))|_{4,A}+
     |F_{}(B_m^*u_x^m(s,\cdot))|_{4,A,1})ds.
$$
 by Lemma \ref{lemmab}. Hence,  
 for $\alpha$ defined as in \eqref{varphi} we have
$$
  K_3 \leq c\int_0^t(t-s)^{-1/2}(|u^m(t,\cdot)|_{4,(-A)^\alpha,1}+
     |(-A)^{-\alpha}D^2u^m(t,\cdot)(-A)^{-1}|_{4,A})ds.
$$
Consequently, taking into account Proposition \ref{E.1}, Proposition \ref{E.2} we find
\begin{eqnarray*}
 K_3 &\leq& 
           c(|\varphi|_{0}+|\varphi|_{2,A,1}+|\Phi|_{0}+|\Phi|_{2,A,1}) \\
     &&\qquad\qquad   \times	    \int_0^t(t-s)^{-1/2}(s^{\alpha-1}+s^{\alpha+\sigma+\varepsilon-3})ds \\
 &\leq& c(|\varphi|_{0}+|\varphi|_{2,A,1}+|\Phi|_{0}+|\Phi|_{2,A,1})(1+t^{\alpha+\sigma+\varepsilon-5/2}).
\end{eqnarray*}
Now the conclusion follows by gathering all the estimates on $K_1$, $K_2$, $K_3$, $K_4$ and by taking into account that $\alpha+\sigma+\varepsilon-5/2 > -1/2$.
\end{proof}
The proof of the following proposition is the same of Proposition 3.6 in \cite{DO};
in addition we have the hamiltonian term, but it is treated by the same arguments of Proposition \eqref{E.1}, \eqref{E.2}, \eqref{E.3}. 
\begin{proposition} \label{E.5}
For all $t$, $s>0$, $m\in\Nset$ and $x\in D(A)$ we have:
\begin{eqnarray*}
|u^m(t,x)-u^m(s,x)| &\leq& c(|\varphi|_{0}+|\varphi|_{\mathcal{E}_{2}}+
      |\Phi|_{0}+|\Phi|_{\mathcal{E}_{2}})(1+|Ax|)^{6} \\
      &&\quad \times \big(|t-s|^{g/2} + |t-s|^{1/2} + |A(\text{\rm e}^{tA}-\text{\rm e}^{sA})x|\big).
\end{eqnarray*}
\end{proposition}
%
%
%
\section{Construction of a solution}
In order to prove existence of solutions of \eqref{R_t} we proceed by a compactness method.
Let $\delta >0$ and set $K_\delta=\{x\in D(A): |Ax|\leq \delta\}$.
\begin{theorem} \label{soluzmild}
Let us assume that \eqref{Q1}, \eqref{Q2}, \eqref{B}, \eqref{varphi} hold\footnote{Existence of solutions for equation \eqref{HJB} may be proved by assuming  that for some $p>0$ the functions $\varphi$, $\Phi$ satisfy $\varphi(x)$, $\Phi(x)\leq c(1+|x|^p)$, $x\in H$, $p>0$. 
However, since  we shall need that in the control problem $\varphi$, $\Phi$ are bounded, we have omitted this case.}. 
Then there exists  a subsequence $\{u^{m_n}(t,x)\}$ of $\{u^m(t,x)\}$ and a continuous function $u:[0,T]\times D(A) \to \Rset$ such that for all $t_0\in(0,T)$, $\delta>0$ the following statements are satisfied:
\begin{itemize}
\item[\rm(i)] $u(t,\cdot)\in C_b(D(A);\Rset)$ for all $t\in [0,T]$. 
Moreover,
\begin{equation} \label{e8.1}
   \lim_{m_n\to \infty} u^{m_n}(t,x) = u(t,x),
\end{equation} uniformly on $[0,T]\times K_\delta$;
%
\item[\rm(ii)] $u(t,\cdot)\in C_2^1(D(A);\Rset)$ for all $t\in(0,T]$. 
Moreover,
\begin{equation} \label{e8.2} 
 \lim_{m_n\to \infty} |(-A)^{-\alpha}\big(Du^{m_n}(t,x) - Du(t,x)\big)| =0,
\end{equation}
uniformly on $[t_0,T]\times K_\delta$;
\item[\rm(iii)] for any $(t,x)\in(0,T]\times D(A)$ there exists the directional derivative             $D_hu(t,x):=Du(t,x)\cdot h$ in any direction $h\in D((-A)^{1/2})$ and
\begin{equation*} 
   \lim_{m_n\to \infty} Du^{m_n}(t,x)\cdot h = Du(t,x)\cdot h,
\end{equation*}
uniformly on $[t_0,T]\times K_\delta$.
Moreover, there exists $c>0$ such that for any $t\in (0,T]$ it holds
\begin{equation*}
    |Du(t,x)\cdot h|\leq ct^{-1/2}(1+|Ax|)^{4}|(-A)^{1/2}h|;
\end{equation*}
\item[\rm(iv)]  
  $u(t,x)$ is a mild solution of the Hamilton-Jacobi-Bellman equation \eqref{HJB}.
\end{itemize}
\end{theorem}
\begin{proof}
Let $\varphi,\Phi \in C_b(D(A);\Rset)\cap \mathcal{E}_{2}(D(A);\Rset)$. 
We deduce from Propositions \ref{E.2}, \ref{E.5} that for any $\delta >0$, $t_0\in (0,T)$, $\sigma \in (3/4,1)$ there exists $C(\delta, t_0, \varphi,g)$ such that $\forall m\in\Nset$, $t,s\geq t_0$ it holds
$$
 |u^m(t,x)-u^m(s,y)|\leq C(\delta, t_0, \varphi,\Phi)\big(|t-s|^{g/2} + |(-A)^{1/2}(x-y)|\big)
$$
and
\begin{equation} \label{convii}
 |(-A)^{-\alpha}\big(Du^m(t,x) - Du^m(t,y)\big)|
    \leq C(\delta, t_0, \varphi,\Phi)|(-A)^{\sigma}(x-y)|.
\end{equation}
Then, by the Ascoli-Arzel\`a Theorem and by a diagonal extraction argument, it follows that there exists a continuous function $u:(0,T]\times D(A)\to \Rset$ such that $u(t,\cdot)\in C_b(D(A);\Rset)$ $\cap$ $C_{2}^1(D(A);\Rset)$ for any $t\in (0,T]$ and a 
subsequence $\{m_n\}_{n\in\Nset}$ of $\Nset$ such that \eqref{e8.1}, \eqref{e8.2} holds, uniformly in $[t_0,T]\times K_\delta$, for any $t_0\in(0,T]$.
Hence (ii) follows.

Let us set $u(0,x)=\varphi(x)$. 
In order to prove that $u:[0,T]\times D(A)\to \Rset$ is continuous, 
it is sufficient to check that \eqref{e8.1} holds uniformly in $[0,T]\times K_\delta$.
Before proving this, we prove (iii).
Let $m,n \in \Nset$, $h\in D((-A)^{1/2})$ and $h' \in D(A)$.
Then for any $t\in (0,T]$ we have
\begin{multline*}
   |(Du^m(t,x)-Du^n(t,x))\cdot h|\\  
   \shoveleft{ \leq|Du^m(t,x)\cdot (h-h')|+|(Du^m(t,x)-Du^n(t,x))\cdot h'|+|Du^n(t,x)\cdot (h-h')|} \\ 
    \shoveleft{\leq t^{-1/2}c(\varphi)(1+|Ax|)^{4}|(-A)^{1/2}(h-h')|}\\
    +|(-A)^{-\alpha}(Du^m(t,x)-Du^n(t,x))||(-A)^{\alpha}h'|,
\end{multline*}
by Proposition \ref{E.3} and \eqref{convii}.
Hence, 
since $D((-A)^{1/2})$ is dense in $D(A)$, the sequence $\{Du^{m_n}(t,x)\cdot h\}_{m\in\Nset}$ is Cauchy in $\Rset$, uniformly in $[t_0,T]\times K_\delta$.
We denote the limit by $Du(t,x)\cdot h$: 
necessarily it coincides with the Gateaux derivative $D_hu(t,x)$ along the direction $h$. Then (iii) is proved.

Now we prove  (iv).
To do this, we shall check that the right-hand side of \eqref{R_tm} converges to the right-hand side of \eqref{R_t}.
By Proposition \eqref{E.3} we have
\begin{equation*}
   |(-A)^{-1/2}Du^m(t,x)|_{4,A}, |(-A)^{-1/2}Du(t,x)|_{4,A}\leq c(\varphi,g)t^{-1/2}.
\end{equation*}
Then by Lemma \ref{lemma2.1} and (iii) we find
\begin{equation*}
  \big(b_m(Z^m(t-s,x)),Du^m(s,Z^m(t-s,x))\big) \leq c(\varphi,g)t^{-1/2}(1+|AZ^m(t-s,x)|)^{5}.
\end{equation*}
Similarly, by \eqref{B} and by Proposition \ref{E.1} we find
\begin{equation*}
F(B_m^*u_x^m(s,Z^m(t-s,x))) \leq s^{-\alpha}c(\varphi,g)(1+|AZ^m(t-s,x)|)^{3}.
\end{equation*}
Then all the integrals in  \eqref{R_t}, \eqref{R_tm} are well defined, and by letting $m_n\to \infty$ in \eqref{R_tm}.
Statement (iv) follows by the dominated convergence theorem and by the well known properties of the Galerkin approximations of the Ornstein-Uhlembeck semigroup.
To complete the proof of (i), it is sufficient to notice that the convergence is uniform in $[0,T]\times K_\delta$.
\end{proof}
\subsection{Martingale solutions of the controlled equation and main result}
The goal of this section is to prove that there exists a solution, in a suitable sense, of equation \eqref{contr}.
We shall prove that the process $X(t,x)$ is a solution of equation \eqref{contr} in the following sense (cf \cite{FG94}): 
\begin{definition}\emph{
 We say there exists a martingale solution of  equation \eqref{contr} if there exists a stochastic basis $(\Omega, \mathcal{F},\{\mathcal{F}_t\}_{t\in[0,T]}, \PP)$, a cylindrical Wiener process $W$ on the space $H$ and a progressively measurable process $X(t,x):[0,T]\times\Omega\to H$ with $\PP$-a.s. paths
$$
  X(\cdot,x)(\omega) \in C([0,T];D((-A)^{-1}))) \cap L^\infty([0,T];H)\cap L^2([0,T];V)
$$
such that the identity
\begin{eqnarray*}
 &&(X(t,x), y) +\int_0^t ( AX(s,x)+b(X(s,x))+Bz(s),y) ds\\
&&\quad=( x, y) +\int_0^t ( Q^{1/2}dW(t),y)
\end{eqnarray*}
holds true for all $t\in [0,T]$, $y\in D(A)$.}
\end{definition}

By It\^o's formula it follows easily (see, for instance, \cite{FG94})
\begin{lemma}
There exists $c>0$ such that for any $x\in H$, $m\in \Nset$, $z\in \mathcal{M}_R$ the following estimate holds:
\begin{equation} \label{bound1}
  \EE\left[\sup_{0\leq t\leq T}|X_m(t)|^2 + \int_0^T\|X_m(s)\|^2ds\right] \leq 
    c(1+|x|^2+\textrm{\em Tr}[Q]),
\end{equation}
\end{lemma}

It is well known (see, for instance, \cite{FG94}), that \eqref{bound1} implies 
that the family of laws 
$\{\mathcal{L}(X^m(\cdot,x))\}_{m\in \Nset}$ is tight in  $L_W^2([0,T];D((-A)^\sigma))$ $ \cap$ $ C([0,T];D((-A)^{-1})$ for any $\sigma<1/2$.
Then, for a fixed $\sigma<1/2$, there exists a probability law $\nu_x$ on $L_W^2([0,T];D((-A)^\sigma))$ $ \cap$ $ C([0,T];D((-A)^{-1}))$ and a subsequence $X^{m_k}(\cdot,x)$ such that $\mathcal{L}(X^{m_k}(\cdot,x)) \to \nu_x$ weakly, i.e 
\begin{eqnarray} 
   \lim_{m_k\to \infty}&&\int_{L_W^2([0,T];D((-A)^\sigma)) \cap C([0,T];D((-A)^{-1}))}
    \psi(\xi)\mathcal{L}(X^{m_k}(\cdot,x))(d\xi) \notag\\ \label{convlaw}
  = &&\int_{L_W^2([0,T];D((-A)^\sigma)) \cap C([0,T];D((-A)^{-1}))}\psi(\xi)\nu_x(d\xi),
\end{eqnarray}
for all bounded and continuous $\psi: L_W^2([0,T];D((-A)^\sigma))$ $ \cap$ $ C([0,T];D((-A)^{-1}))$ $\to\Rset$.
Moreover, by the Skorohod theorem (see, for instance, \cite{DPZ}), there exists a probability space $(\Omega_x,\mathcal{F}_x,\PP_x)$ and a stochastic process $X(\cdot,x)$ with law $\nu_x$ such that $X^{m_k}(t,x) \to X(t,x)$  $\PP_x$-a.s..

So, we have proved the following
\begin{theorem} \label{solutioncontr}
For any $z\in \mathcal{M}_R$ and for any $x\in D(A)$ there exists a probability space $(\Omega_x, \mathcal{F}_x,\PP_x)$, a process $X(\cdot,x)$ and a sequence of stochastic processes $\{X^{m_n}(\cdot,x)\}_{k\in\Nset}$ defined on $(\Omega_x, \mathcal{F}_x,\PP_x)$ such that     
\begin{itemize}
\item[\rm (i)] The processes $X^{m_n}(\cdot,x)$ are solutions of the Galerkin approximated equations               \eqref{contrm};
\item[\rm (ii)] The sequence $X^{m_n}(\cdot,x)$ converges $\PP_x$-a.s. to the process $X(\cdot,x)$;
\item[\rm (iii)] $X(\cdot,x)$  is a martingale solution of \eqref{contr}.
\end{itemize}
\end{theorem}
Let us assume that there exists a solution $u(t,x)$ of \eqref{R_t} such that for a subsequence $u^{m_n}(t,x)$ of solutions of \eqref{R_tm} the statements (i)-(iv) of Theorem \ref{solutioncontr} are satisfied and let us denote by $\widetilde{X}^{m_n}(\cdot,x)$ the solution of the $m_n$-Galerkin approximated equation \eqref{loopm}.
Since $P_mz^* \in \mathcal{M}_R$ we can argue as in Theorem \eqref{solutioncontr} to find a subsequence $\{m_n'\}_{n\in\Nset}$ of $\{m_n\}_{n\in\Nset}$ and a probability space $(\Omega_x, \mathcal{F}_x,\PP_x)$ such that the solutions $\widetilde{X}^{m_n'}(\cdot,x)$
converges  $\PP_x$-a.s. to a process $\widetilde{X}(\cdot,x)$.  
In order to prove that $\widetilde{X}(\cdot,x)$ is  a martingale solution of the closed loop equation \eqref{loop} we have to define the optimal control $z^*$ as in \eqref{z^*}. 
Then, since $u(t,x)$ is defined for $x\in D(A)$, we need that the martingale solution $X(t,x)$ of \eqref{contr} belongs to $D(A)$.
To do this, we have the next
\begin{lemma}
For any $\delta \in (1/2,\min\{1+g,1+2\gamma\}]$ there exists a constant $c(\delta)>0$ such that for any $x\in H$, $m\in \Nset$ and $t\in [0,T]$ it holds
$$
  \EE \bigg[\int_0^T \frac{|(-A)^{(1+\delta)/2}X^m(s,x)|^2}{(1+|(-A)^{\delta/2}X^m(s,x)|^2)^{\theta_\delta}}ds\bigg]
    \leq c(\delta)
$$
$$
\text{where %
 $\theta_\delta=\frac{2\delta+1}{2\delta-1}$} 
$$
\end{lemma}
\begin{proof}
We apply the It\^o formula to 
$$
    F_\delta(x)= - \frac{1}{(1+|(-A)^{\delta/2}x|^2)^{\theta_\delta-1}}.
$$
We obtain
\begin{eqnarray*}
  &&\frac{1}{(1+|(-A)^{\delta/2}x|^2)^{\theta_\delta-1}}+
    2(\theta_\delta -1)  \EE\bigg[\int_0^T
    \frac{|(-A)^{(1+\delta)/2}X^m(s,x)|^2}
         {(1+|(-A)^{\delta/2}X^m(s,x)|^2)^{\theta_\delta}}ds\bigg] \\
&&\quad= 2(\theta_\delta -1) \EE\bigg[\int_0^T\frac{(b_m(X^m(s,x)),(-A)^\delta X^m(s,x))}
         {(1+|(-A)^{\delta/2}X^m(s,x)|^2)^{\theta_\delta}}ds\bigg] \\
 &&\qquad+ (\theta_\delta -1) \EE\bigg[\int_0^T
    \frac{((-A)^\delta X^m(s,x),Bz(t))}
         {(1+|(-A)^{\delta/2}X^m(s,x)|^2)^{\theta_\delta}}ds\bigg] \\
 &&\qquad+ (\theta_\delta -1) 
   \EE\bigg[\int_0^T \frac{1}{(1+|(-A)^{\delta/2}X^m(s,x)|^2)^{\theta_\delta}}ds\bigg]
      \text{ Tr}[Q_m(-A)^\delta]\\
  &&\qquad-\theta_\delta(\theta_\delta -1)
  \EE\bigg[\int_0^T \frac{|Q_m^{1/2}(-A)^{\delta/2}X^m(s,x)|^2}
                   {(1+|(-A)^{\delta/2}X^m(s,x)|^2)^{\theta_\delta+1}}ds\bigg]\\
&&\qquad+\EE\bigg[\frac{1}{(1+|(-A)^{\delta/2}X^m(T,x)|^2)^{\theta_\delta-1}}\bigg].
\end{eqnarray*}
By Lemma \ref{lemma2.1} and by interpolation we find
\begin{eqnarray*}
(b_m(X_m),(-A)^\delta X^m) &\leq&
 c|(-A)^{\delta/2}X^m|^{1/2+\delta} |(-A)^{(1+\delta)/2}X^m|^{5/2-\delta} \\
  &\leq& c|(-A)^{\delta/2}X^m|^{2\frac{2\delta+1}{2\delta-1}} + 
    \frac{1}{4} |(-A)^{(1+\delta)/2}X^m|^2.
\end{eqnarray*}
Moreover, since $\delta \leq 1+2\gamma$ and $z$ is bounded we have
\begin{eqnarray*}
  ((-A)^\delta X^m, Bz)&\leq& c|B^*(-A)^\delta X^m|  \leq c|(-A)^{\delta-\gamma}X^m|\\
  &\leq& c |(-A)^{(1+\delta)/2}X^m| \leq c+\frac{1}{4}|(-A)^{(1+\delta)/2}X^m|^2.
\end{eqnarray*}
By \eqref{Q1} we know that $Q_m^{1/2}(-A)^{\delta/2}$ is a bounded operator. 
It follows that the last three terms in the right-hand side are bounded.
We deduce that
\begin{eqnarray*}
  &&\EE\bigg[\int_0^T
    \frac{|(-A)^{(1+\delta)/2}X^m(s,x)|^2}
         {(1+|(-A)^{\delta/2}X^m(s,x)|^2)^{\theta_\delta}}ds\bigg] \\
  &&\quad\leq \EE\bigg[\int_0^T
    \frac{|(-A)^{\delta/2}X^m(s,x)|^{2\frac{2\delta+1}{2\delta-1}}}
         {(1+|(-A)^{\delta/2}X^m(s,x)|^2)^{\theta_\delta}}ds\bigg] + c(\delta).
\end{eqnarray*}
Finally, since $\theta_\gamma= \frac{2\delta+1}{2\delta-1}$, the result follows.
\end{proof}
Arguing as in Lemma 7.5 of \cite{DPD03} it follows the next
\begin{lemma} \label{solutioncontrA}
Under the assumptions of Theorem \eqref{solutioncontr} the sequence of processes $\{X^{m_n}(\cdot,x)\}_{n\in\Nset}$, $X(\cdot,x)$ also satisfies 
\begin{equation} \label{convD(A)}
   \lim_{m_n\to\infty}X^{m_n}(t,x) =X(t,x)\text{ in } D(A),\, dt\times \PP_x\text{-a.s in } [0,T]\times\Omega_x.
\end{equation}
\end{lemma}

\subsubsection*{Main result}
For any sequence $\{m_n\}_{n\in\Nset}\subset \Nset$ and any $z\in \mathcal{M}_R$ we denote by $MS(\{m_n\},z)$ the set of all martingale solutions $X(\cdot,x)$ of \eqref{contr} such that for some subsequence
$\{m_n'\}_{n\in\Nset}\subset\{m_n\}_{n\in\Nset}$ we have
$$
  X(\cdot,x) = \lim_{m_n'\to \infty} X^{m_n'}(\cdot,x)\quad \text{in law},
$$
where $ X^{m_n'}(\cdot,x) $ is the solution of the the $m_n'$-Galerkin approximated equation \eqref{contrm}.
We have the next
\begin{theorem}
Let us assume that the functions $u(t,x)$, $\{u^{m_n}(t,x)\}_{n\in\Nset}$ satisfy statements (i)-(iv) of Theorem \ref{soluzmild}. 
Then, for any $x\in D(A)$ there exists a martingale solution $\widetilde{X}(\cdot,x)$ of the closed loop equation \eqref{loop}. 
Moreover, the control $z^*$ defined by
$$
   z^*=-D_pF(B^*u_x(T-t,\widetilde{X}(t,x)))
$$
verifies
$$
  u(T,x) =J(\widetilde{X}(\cdot,x),z^*) = 
  \min_{\substack{ z\in{\cal M}_R\\X(\cdot,x)\in MS(\{m_n\},z)}}
    J(X(\cdot,x),z).
$$

\end{theorem}
\begin{proof}
Arguing as in Theorem \ref{solutioncontr}, there exists a subsequence $\{m_n'\}_{n\in\Nset}$ of $\{m_n\}_{n\in\Nset}$, a probability space $(\widetilde{\Omega}_x, \widetilde{\mathcal{F}}_x,\widetilde{\PP}_x)$ and a sequence of processes $\{\widetilde{X}^{m_n'}(\cdot,x)\}_{n\in\Nset}$, $\widetilde{X}(\cdot,x)$ defined on $\widetilde{\Omega}_x$ which satisfies statements (i), (ii) of Theorem \ref{solutioncontr} and \eqref{convlaw}.
Moreover, by Lemma \ref{solutioncontrA} we can assume that \eqref{convD(A)} holds.
Consequently, by   statement (ii) of Theorem \ref{soluzmild} it is easy to see that
$$
  \lim_{m_n'\to \infty}B_{m_n'}^*u_x^{m_n'}(T-t,\widetilde{X}^{m_n'}(t,x)) =  B^*u_x(T-t,\widetilde{X}(t,x)),\quad dt\times\widetilde{\PP}_x\text{-a.s.}.
$$
Hence, since $|D_pF| \leq R$, we can apply the dominated convergence theorem in order to have
\begin{eqnarray}
 && \lim_{m_n'\to \infty}\int_0^t (D_pF(B_{m_n'}^*u^{m_n'}(T-s,\widetilde{X}^{m_n'}(s,x))),y) ds  \notag \\
 &&\quad = \int_0^t( D_pF(B^*u(T-s,\widetilde{X}(s,x))),y) ds, \label{convzz}
\end{eqnarray}
$\widetilde{\PP}_x\text{-a.s.}$, for all $t\in [0,T]$ and all $y\in D((-A)^{-1})$.
Now, arguing as in Theorem \ref{solutioncontr}, it follows that $\widetilde{X}(\cdot,x)$ is a martingale solution of \eqref{loop}.

For the second statement we notice that it holds 
$$
  u(T,x)-\widetilde{\EE}_x\big[\varphi(\widetilde{X}(T,x))\big] = \lim_{m_n\to\infty}\Big( u^{m_n'}(T,x)-\widetilde{\EE}_x\big[\varphi(\widetilde{X}^{m_n'}(T,x))\big]\Big).
$$ 
Since $u^m(t,x)$ is solution of \eqref{HJBm}, by a standard computation based on the It\^o formula (cf \cite{DPZ02}) we have, for any $m\in \Nset$,
$$
 u^m(T,x)=\widetilde{\EE}_x\bigg[\int_0^T\big(\Phi(\widetilde{X}^m(s,x))+\frac{1}{2}|z^*_m(s)|^2\big) ds +\varphi(\widetilde{X}^m(T,x))\bigg].
$$
Then, taking into account \eqref{convzz} and that $\varphi$, $\Phi$ are bounded, 
we can apply the dominated convergence Theorem to yield 
$$
  u(T,x) = J(z^*)(x).
$$
Now notice that for any $n\in\Nset$ the process $\widetilde{X}^{m_n'}(\cdot,x)$ is the optimal state for the $m_n'$-Galerkin approximated control problem \eqref{contrm}, which implies 
$$
 u^{m_n'}(T,x)\leq 
 \EE_x\bigg[\int_0^T\big(\Phi(X^{m_n'}(s,x))+\frac{1}{2}|P_{m_n'}z(s)|^2\big)ds+\varphi(X^{m_n'}(s,x))\bigg],
$$
for all $z\in \mathcal{M}_R$.
Then, by letting $m_n'\to \infty$, we conclude the proof.
\end{proof}
{\bf Acknowledgements}
I would like to thank A.Debussche for introducing me to this beautiful subject.

\end{document}